\documentclass[12pt]{amsart}

\usepackage{color,epsfig,amsfonts,fancyhdr,ifthen}
\usepackage[letterpaper,left=1.25in,right=1.25in,top=1.25in,bottom=1.
25in]{geometry}

\newtheorem{thm}{Theorem}[section]
\newtheorem*{thm*}{Theorem}
\newtheorem*{prop*}{Proposition}
\newtheorem{prop}[thm]{Proposition}
\newtheorem{cor}[thm]{Corollary}
\newtheorem{lem}[thm]{Lemma}

\newtheorem{conj}[thm]{Conjecture}
\newtheorem{defn}[thm]{Definition}
\newtheorem{rmk}{Remark}
\newtheorem{example}{Example}
\newtheorem{question}{Question}

\newcommand{\BP}{Blaschke product }
\newcommand{\BPNSP}{Blaschke product}
\newcommand{\BPS}{Blaschke products }
\newcommand{\BPSNSP}{Blaschke products}

\newcommand{\PP}{\ensuremath{\mathbb{P}}}
\newcommand{\C}{\ensuremath{\mathbb{C}}}
\newcommand{\R}{\ensuremath{\mathbb{R}}}
\newcommand{\T}{\ensuremath{\mathbb{T}}}
\newcommand{\D}{\ensuremath{\mathbb{D}}}
\newcommand{\N}{\ensuremath{\mathbb{N}}}
\newcommand{\B}{\ensuremath{\mathcal B}}
\newcommand{\KK}{\ensuremath{\mathcal K}}

\newcommand{\ra}{\ensuremath{\rightarrow}}
\newcommand{\dt}{\ensuremath{d_{\rm top}}}
\newcommand{\da}{\ensuremath{d_{\rm alg}}}
\newcommand{\mutor}{\ensuremath{\mu_{\rm tor}}}
\newcommand{\tr}{{\rm tr}}
\newcommand{\supp}{{\rm supp}}
\newcommand{\sm}{\setminus}

\newcommand{\im}{{\rm Im}}
\newcommand{\dist}{{\rm dist}}
\newcommand{\HH}{\mathcal{H}}

\newcommand{\loc}{{\rm loc}}

\begin{document}

\author{Enrique R. Pujals and Roland K. W. Roeder}

\title{Two-dimensional Blaschke products:  degree growth and ergodic consequences}


\date{\today}

\begin{abstract}
We study the dynamics of Blaschke products in two dimensions, particularly the 
rates of growth for the degrees of iterates and the corresponding implications for the ergodic
properties of the map.
\end{abstract}

\maketitle
\markboth{\textsc{E. R. Pujals and R.K.W. Roeder}}
  {\textit{Two dimensional Blaschke products}}

\section{Introduction}

For dominant rational maps of compact, complex, Kahler manifolds there is a
conjecture specifying the expected ergodic properties of the map depending on
the relationship between the rates of growth for certain degrees under
iteration of the map.  (See Conjecture \ref{CONJ:ENTROPY}, as presented in
\cite{GUEDJ_ERGODIC}, as well as the results towards this conjecture
\cite{GUEDJ_LARGE_DEGREE,DDG1,DDG2,DDG3}.) We observe that the two-dimensional
\BPS fit naturally within this conjecture, having examples from each of the
three cases that the conjecture gives for maps of a surface.    We then consider
the dynamics of \BPS from these dramatically distinct cases, relating it to the
behavior predicted by this conjecture.

Furthermore, generic (in an appropriate sense) \BPS do not have the  Julia set
contained within $\T^2$.  Rather, ``the majority of it'' is away from $\T^2$
within the support of the measure of maximal entropy.  This is very different
from the case of $1$-dimensional \BPS for which the Julia set is the unit
circle (see below).

A (finite) Blaschke product is a map of the form
\begin{eqnarray}\label{EQN:BP1}
E(z) = \theta_0 \prod_{i=1}^n \frac{z-e_i}{1-z \overline{e_i}},
\end{eqnarray}
\noindent
where $n \geq 2$, $e_i \in \C$ for each $i=1,\ldots,n$, and $\theta_0 \in \C$ with $|\theta_0| = 1$.
The simplest dynamical situation occurs if one restricts that $|e_i| < 1$ for $i=1\ldots n$.  It implies that the Julia set 
$J_E$ is contained within the unit circle $\T^1$. 

In this paper we study \BPS  in two variables generalizing this situation.  Let
\begin{eqnarray}\label{EQN:BP2}
f(z,w) = \left(\theta_1 \prod_{i=1}^{m} \frac{z-a_i}{1-\bar a_i z} \prod_{i=1}^{n} \frac{w-b_i}{1-\bar b_i w},
\theta_2 \prod_{i=1}^{p} \frac{z-c_i}{1-\bar c_i z} \prod_{i=1}^{q} \frac{w-d_i}{1-\bar d_i w}\right),
\end{eqnarray}
\noindent
with $|\theta_1|=|\theta_2| = 1$ and all of the zeros $a_1,\ldots,d_q$ of modulus less than one.
We will often denote the corresponding $1$-variable \BPS by $A(z), B(w), C(z),$
and $D(w)$.  Such maps were introduced in \cite{BP_PS}.

Note that if one allows some of the zeros $e_i$ of a one-variable \BP
(\ref{EQN:BP1}) to have modulus greater than $1$, a much more complicated
structure for the Julia set can occur \cite{ROE}.  We do not consider the
$2$-variable analog of that situation in this paper, but we expect that it may be interesting for
further study.

We describe the degrees of a given \BP $f$
by a matrix
\begin{eqnarray*}
N = \left[\begin{array}{cc} m & n \\ p & q \end{array} \right].
\end{eqnarray*}
\noindent
(As in \cite{BP_PS}, we assume that $m,n,p,$ and $q$ are greater than or equal
to $1$).

Given any matrix of degrees $N$, any choice of rotations
$\theta_1,\theta_2$, and any zeros $a_1,\ldots,a_m$, $b_1,\ldots,b_n$,
$c_1,\ldots,c_p$, and $d_1,\ldots,d_q$ (all of modulus less than $1$) there is
a \BPNSP.  We denote the space of all such \BPS by $\B_N$ and we will call any $f
\in \B_N$ a {\em Blaschke Product associated to $N$}.  We typically will use
the notation $\sigma \in \D^{m+n+p+q}$ to represent the collection of zeros
$a_1,\ldots,d_q$.  Notice that $\B_N$ can be identified with $\D^{m+n+p+q} \times \T^2$,
an identification that we use when discussing sets of full measure on $\B_N$.

In the case that all of the zeros are equal to  $0$, a
2-dimensional \BP becomes a {\em monomial map}
\begin{eqnarray}
f(z,w) = (z^m w^n, z^p w^q),
\end{eqnarray}
\noindent
whose dynamics  was studied extensively in \cite{FAVRE_MON,HP_MON}.  For any $N$ we will
call this map the {\em monomial map associated to $N$}.  (It is also interesting to
note that monomial maps occur frequently ``outside of dynamical systems'', for
example in the description of cusps for Inoue-Hirzebruch surfaces \cite{DLOU}).

\vspace{0.05in}
One nice reason to study \BPS is that they preserve the unit torus $\T^2 :=
\{(z,w) \,\, : \,\, |z|=|w| = 1\}$.  The monomial map associated to $N$ induces
a linear map on $\T^2$. If $\det N > 0$, this is an orientation preserving
local diffeomorphism of topological degree $\det N$.  (The topological degree
is the number of preimages of a generic point for $f_{|\T^2}$).  {\em
Throughout the paper we will assume $\det N > 0$.} Furthermore, the action on
$\pi_1(\T^2)$ is described by $N$, in terms of the obvious choice of
generators.

Any $f \in \B_N$ is homotopic on $\T^2$ to this monomial map and therefore has
the same action on $\pi_1(\T^2)$ and the same topological degree.  However, it
may fail to be a local diffeomorphism.  

We will often consider the special case in which $f_{|\T^2}$ is an orientation
preserving diffeomorphism of $\T^2$.   We call such an $f$ as a {\em \BP
diffeomorphism} (although generally it is only a diffeomorphism on
$\T^2$, not globally on $\PP^2$).  \BP diffeomorphisms were studied extensively
in \cite{BP_PS} and they can only occur if $\det N = 1$.  The
corresponding monomial map induces a linear Anosov map on $\T^2$ and a \BP
whose zeros are sufficiently small will also be a \BP diffeomorphism, inducing
an Anosov map on $\T^2$.  For any \BP diffeomorphism,  the
restriction $f|_{\T^2}$ has topological entropy $\log(c_+(N))$, where
$c_+(N)$ is the largest eigenvalue of $N$. There is also a unique invariant measure $\mutor$ of
maximal entropy for $f_{|\T^2}$.  (See Appendix \ref{APP:ENTROPY_T2}.)

\vspace{0.05in}

A rational map $g:\PP^2 \ra \PP^2$ can be lifted to a system of three homogeneous equations on $\C^3$ having no common factors.  The {\em algebraic degree} $\da(g)$ is
the common degree of these homogeneous equations.  In some
cases, the degree of iterates drops, $\da(g^n) < \left(\da(g)\right)^n$, because a
common factor appears in the homogeneous equations for $g^n$. (See \cite{FS}).  However, a
limiting degree called {\em the first dynamical degree}
\begin{eqnarray}
\label{EQN:DYN_DEG1}
\lambda_1(g) = \lim_{n\rightarrow \infty} (\da(g^n))^{1/n}
\end{eqnarray}
always exists, describing the asymptotic rate of growth in the sequence
$\{\da(g^n)\}$, \cite{RUSS_SHIFF}.  Note that $\lambda_1(g) \leq \da(g)$.
In \S \ref{SEC:HOW_TO_COMPUTE_DYN_DEG} we briefly describe a common technique for computing $\lambda_1(f)$.

The ergodic properties of $g$ are believed, see \cite[Conj. 3.2]{GUEDJ_ERGODIC}, to
depend heavily on the relationship between $\lambda_1(g)$ and the topological
degree $\dt(g)$. (Here, $\dt(g)$ is defined as the
number of preimages under $g$ of a generic point from $\PP^2$.) Actually, the
conjecture stated in \cite{GUEDJ_ERGODIC} is far more general, pertaining to
dominant maps of Kahler manifolds $X$ of arbitrary dimension.  We provide a
brief summary in the case that $X$ is a surface.

\begin{conj}\label{CONJ:ENTROPY}
The ergodic properties of a rational map $g:X \rightarrow X$ are believed to fall into three cases:\\
\begin{itemize}
\item {\bf Case I:  Large topological degree:} $\dt(g) > \lambda_1(g)$.
 This case has been solved by \cite{GUEDJ_LARGE_DEGREE} (see also \cite{RUSS_SHIFF}) where it was shown
that there is an  ergodic invariant measure $\mu$ of maximal
entropy $\log(\dt(f))$. The measure $\mu$ is not supported on hypersurfaces, it
does not charge the points of indeterminacy, and the repelling points of $f$ are
equidistributed according to this measure.  It is the unique measure of maximal entropy.

\item {\bf Case II: Small topological degree:} $\dt(g) < \lambda_1(g)$. 
It is believed that there is an ergodic invariant measure $\mu$ of maximal
entropy $\log(\lambda_1(g))$ that is not supported on hypersurfaces and does
not charge the points of indeterminacy.  Saddle-type points are believed to be
equidistributed according to this measure.  It is the unique measure of maximal entropy.

A recent series of preprints
\cite{DDG1,DDG2,DDG3} has appeared where it is proven that these expected
properties (except for uniqueness of $\mu$) hold, provided that certain technical hypotheses are met. 

\item {\bf Case III: Equal degrees:} $\dt(g) = \lambda_1(g)$.  Little is known or conjectured in this case.

\end{itemize}
\end{conj}

\vspace{0.1in}
\begin{rmk}\label{RMK:MONOMIALS}
Suppose that $f$ is the monomial map associated to $N$.  According to \cite{FAVRE_MON}, $\lambda_1(f) = c_+(N)$ and $\dt(f) = \det N$.  Therefore,
by choosing $N$ appropriately we can find monomial maps in each of the three cases from Conjecture \ref{CONJ:ENTROPY}.
\end{rmk}

We now summarize the main results of this paper:

In \S \ref{SEC:DEGREES} we prove

\begin{thm}
\label{THM:DYN_DEGREE}
Any $f \in \B_N$  has the same dynamical degree as the monomial map associated to $N$.  That is:
\begin{eqnarray*}
\lambda_1(f) = c_+(N) =  \frac{m+q+\sqrt{(m-q)^2+4np}}{2},
\end{eqnarray*}
\noindent
where $c_+(N)$  is the leading eigenvalue of $N$.
\end{thm}

In \S \ref{SEC:LARGE_DEGREE} we consider \BPS falling into Case I of 
Conjecture \ref{CONJ:ENTROPY}.

Notice that $\dt(f) \geq \dt(f_{|\T^2})
= \det N$, so that $\dt(f)$ is greater than or equal to the topological degree of the monomial map associated to $N$.
In particular, if the monomial map associated to $N$ falls into Case I of the conjecture (i.e. $\det N > c_+(N)$) then so does
every other $f \in \B_N$.

On the other hand, for any $N$, most \BPS fall into Case I:

\begin{thm}\label{THM:GENERICALLY_LARGE_DEGREE}
For any matrix of degrees $N$
there is an open dense set of full measure $\hat \B_N \subset \B_N$ so that if
$f \in \hat \B_N$ then $\dt(f) = mq+np > \lambda_1(f)$.
\end{thm}
\noindent

The results from \cite{GUEDJ_LARGE_DEGREE} apply, giving the existence
of a unique measure of maximal entropy $\mu$ having entropy $\log(\dt(f))$.
For particular choices of $f$
we can have $\supp(\mu) \subset \T^2$.  However,
if $f \in \hat \B_N$, this measure does not charge the invariant torus
$\T^2$.  Furthermore, in certain situations, an analysis of the dynamics near $\T^2$ allows one to see
that $\supp(\mu)$ is isolated away from $\T^2$.  

In \S \ref{SEC:SMALL_DEGREE} we consider \BPS falling into Case
II of Conjecture \ref{CONJ:ENTROPY}.  If $\det N < c_+(N)$, this occurs for
the monomial maps associated to $N$, as well as certain non-generic $f \in
\B_N$.

Many of the examples in \S \ref{SEC:SMALL_DEGREE} induce a diffeomorphism of
$\T^2$ in which case there is an invariant measure $\mutor$ supported on $\T^2$
of entropy $\log c_+(N) = \log \lambda_1(f)$.  As a consequence of the bound on
entropy provided in \cite{DS_ENTROPY}, we find

\begin{prop}\label{PROP:MAX_ENTROPY_IN_T2}
Let $f$ be a \BP diffeomorphism of small topological degree $\dt(f) < \lambda_1(f)$.  Then, $f:\PP^2 \rightarrow \PP^2$ has a measure of maximal entropy $\mutor$ supported within $\T^2$.
\end{prop}

We do not know if $\mutor$ is the unique measure of maximal entropy in all of $\PP^2$ for these
\BP diffeomorphisms.  Furthermore, it would also be
interesting to see how these maps fit within the framework presented in
\cite{DDG1,DDG2,DDG3}.

In \S \ref{SEC:EQUAL_DEGREES} we briefly consider the case of \BPS
falling into Case III of Conjecture \ref{CONJ:ENTROPY}.

We conclude with Appendix \ref{APP:ENTROPY_T2} by proving basic facts about the entropy of \BP diffeomorphisms in $\T^2$.

\subsection*{Acknowledgments}

We thank the referee for many
helpful suggestions, including ideas that allowed us to dramatically improve
the results stated in Theorem \ref{THM:DYN_DEGREE} and Proposition \ref{PROP:PREIMAGES_OF_T2}.

We thank Jeffrey Diller for introducing the first author to the technique for
computing dynamical degrees that is described in \S
\ref{SEC:HOW_TO_COMPUTE_DYN_DEG} and for other helpful discussions.  John Hubbard informed the
authors of the relationship between monomial maps and Inoue-Hirzebruch surfaces
and he provided many other helpful suggestions. In addition, we have benefited from
interesting conversations with Eric Bedford, Tien-Cuong Dinh, Romain Dujardin,
Michael Shub, and Nessim Sibony.

\section{A standard technique for computing dynamical degree}\label{SEC:HOW_TO_COMPUTE_DYN_DEG}
This section describes the work of many other authors (see the references
within) and none of it is original to this paper.  We provide it as a brief
summary of the technique that we will use for computing the dynamical degree of
\BPSNSP.

One can recasts the dynamical
degree (\ref{EQN:DYN_DEG1}) as:
\begin{eqnarray}\label{EQN:DYN_DEG2}
\lambda_1(f) = \limsup(r_1((f^n)^*))^{1/n},
\end{eqnarray}
\noindent
where $r_1((f^n)^*)$ is the spectral radius of the linear action of $(f^n)^*$ on
$H_a^{1,1}(X,\mathbb{R})$. Here, $H_a^{1,1}(X,\mathbb{R})$ is the part of the $(1,1)$ cohomology that is
generated by algebraic curves in $X$, see \cite[Prop 1.2(iii)]{GUEDJ_LARGE_DEGREE}. (The
cohomology class $[D]$ of an algebraic curve $D$ is taken in the sense of closed-positive $(1,1)$ currents.)
When $X = \PP^2$ this definition agrees with (\ref{EQN:DYN_DEG1}) and
this new definition is invariant under birational conjugacy (see \cite[Prop
1.5]{GUEDJ_LARGE_DEGREE}).

\begin{defn}
A rational mapping $f: X \rightarrow X$ of a Kahler surface $X$ is
called {\em algebraically stable} if there is no integer $n$ and no hypersurface
$V$ so that each component of $f^n(V)$ is contained within the indeterminacy set $I_f$.
\end{defn}

\noindent
For the case $X = \PP^2$, see \cite[p. 109]{S_PANORAME} and more generally, see
\cite{DF_BIRATIONAL}.

If
$f:X \ra X$ is algebraically stable then, according to \cite[Thm
1.14]{DF_BIRATIONAL}, one has that the action of $f^*:H_a^{1,1}(X,\mathbb{R})
\ra H_a^{1,1}(X,\mathbb{R})$ is well-behaved: $(f^n)^* = (f^*)^n$.  In this
case, (\ref{EQN:DYN_DEG2}) simplifies to 
\begin{equation}\label{EQN:DEGREES_AGREE}
\lambda_1(f) =r_1(f^*).
\end{equation}
\noindent
Therefore, {\em computation of dynamical degree for an algebraically stable mapping reduces to the study of a single iterate.}

If $f:\PP^2 \ra \PP^2$ that is not algebraically stable, a
typical way to compute $\lambda_1(f)$ is as follows.  One tries to find an
appropriate finite sequence of blow-ups at certain points in $\PP^2$ in an
attempt to obtain a new surface $X$ on which the extension $\tilde{f}$ of
$f$ is algebraically stable.  Note that in this approach $\tilde{f}$ and $f$ are birationally conjugate
using the canonical projection $\pi: X \ra \PP^2$, and hence $\lambda_1(f) = \lambda_1(\tilde{f})$.

A surface $X$ that is birationally equivalent to $\PP^2$ is called a {\em
rational surface}.  In this paper we will always construct $X$ using the
strategy described in the previous paragraph, so it will be an ongoing
assumption that any surface $X$ is rational (unless otherwise explicitly
stated).  In this case, $H^{1,1}_a(X)$ coincides with the full
cohomology $H^{1,1}(X)$, allowing us a further simplification.

Suppose that one has created such a new surface $X$ so that
$\tilde{f} : X \ra X$ is algebraically stable.  Then, 
$\lambda_1(f) = \lambda_1(\tilde{f}) = r_1(\tilde{f}^*)$, where
$r_1(\tilde{f}^*)$ is  the spectral radius of the action
$\tilde{f}^*:H^{1,1}(X,\mathbb{R}) \ra H^{1,1}(X,\mathbb{R})$.    This latter
number can be computed by considering the pull-backs $f^*$ of an appropriate
finite set of curves that form a basis for $H^{1,1}(X,\mathbb{R})$.  Nice
descriptions of this procedure and explicit examples are demonstrated in
\cite{BEDFORD_KIM_JGA,BEDFORD_DILLER,BEDFORD_KIM} and the references therein.
(The latter two of these references work in terms of $Pic(X)$, rather
than $H^{1,1}(X)$, but the technique is essentially the same.)

In fact, such a modification does not exist for all rational maps.  In
\cite{FAVRE_MON} it was  shown that for certain monomial maps (with some
negative powers) there is no finite sequence of blow-ups that one can do,
starting with $\PP^2$, in order to obtain a surface $X$ on which the map is
algebraically stable.  However, in the case that $f$ is a monomial map with all
positive powers (as assumed in this paper) it was shown in \cite{FAVRE_MON}
that one can always find a toric surface $\check{X}$ on which $f$ becomes
algebraically stable.  In this case, $\check{X}$ is obtained first by
blowing-up $\PP^2$ {\em and then extending to a ramified cover} (so that it is
typically no longer a rational surface). 
See Question \ref{QUESTION:STABLE_MODEL} at the end of \S \ref{SUBSEC:EQUALITY_DENSE_OPEN_SET}.

\section{Computation of dynamical degree for \BPS} \label{SEC:DEGREES}

In this section we prove Theorem \ref{THM:DYN_DEGREE}, which states that for any \BP $f \in \B_N$ we have $\lambda_1(f) = c_+(N)$.

We employ the following strategy: In \S \ref{SUBSEC:LOWER_BOUND} we obtain a
lower bound $\lambda_1(f) \geq c_+(N)$ for all $f \in \B_N$.  It will be a
consequence of the dynamics of $f|_{\T^2}$.  Then, in \S
\ref{SUBSEC:EQUALITY_DENSE_OPEN_SET} we use the strategy described in \S
\ref{SEC:HOW_TO_COMPUTE_DYN_DEG} to find a dense set of full measure $\B'_N
\subset \B_N$ on which $\lambda_1(f) = c_+(N).$  In \S
\ref{SUBSEC:EQUALITY_EVERYWHERE} we combine the results of \S
\ref{SUBSEC:LOWER_BOUND} and \S \ref{SUBSEC:EQUALITY_DENSE_OPEN_SET} to show
$\lambda_1(f) = c_+(N)$ everywhere.

\subsection{Lower bound}\label{SUBSEC:LOWER_BOUND}

\begin{prop}
\label{PROP:DYN_DEG_LOWER}
For any \BP $f \in \B_N$ we have that $\lambda_1(f) \geq c_+(N)$.
\end{prop}

\begin{proof}
We use the definition given in Equation (\ref{EQN:DYN_DEG1}) for $\lambda_1(f)$.

Consider the basis $\{[\gamma_1],[\gamma_2]\}$ for $H_1(\T^2)$ generated by the unit circle
$\gamma_1$ in the plane $w=0$ and the unit circle $\gamma_2$ in the plane $z=0$.
As noted earlier, the action $f_*:H_1(\T^2) \ra H_1(\T^2)$ with respect to this basis is given
by multiplication by the matrix $N$.

We will show that $\da(f^n) \geq ||N^n||_\infty$, i.e.  that $\da(f^n)$ grows
at least as fast as the largest element of $N^n$.  This suffices to prove the
assertion since $||N^n||_\infty \geq a \cdot c_+(N)^n$ for some positive constant $a$.

Notice that $f_*$ acts ``stably'' on $H_1(\T^2)$ in the sense that the action of $f^n_*$ is given by
$N^n$ with respect to the previously mentioned basis.  Consider now the largest element of $N^n$,
which we suppose (for the moment) is the $(1,1)$ element.  Then $f^n_*([\gamma_1]) = k
[\gamma_1]$ where $k \geq a \cdot c_+(N)^n$.

Write $f$ in affine coordinates $(z,w)$.
We will show that the first coordinate of $f^n$ is a rational function of degree at least $k$ in
$z$.  This is sufficient to give that any homogeneous expression for $f^n$ has degree at least
$k$, as well.

Let $\pi$ be the projection $\pi(z,w) = z$ so that the first coordinate of
$f^n$ is given by $\pi \circ f^n$.  Also let $\iota(z) = (z,1)$.  The iterate
$f^n$ is holomorphic on the open bidisc $\D \times \D$ because $f^n$ forms a
normal family there. Then $\pi \circ f^n \circ \iota: \overline \D \ra
\overline \D$ is a holomorphic function preserving the unit circle.  By the
previous homological considerations this map has degree $k$ on the circle, so
by a standard theorem it must be a (one variable) \BP of degree $k$ with no
poles inside of $\D$.
This gives a lower bound for the degree in $z$ of the first coordinate of $f$
by $k \geq a \cdot c_+(N)^n$.

In the case that some other element than the $(1,1)$ element of $N^n$ were largest,
an identical proof works by choosing $\iota$ to be the appropriate inclusion and $\pi$ to be the
appropriate projection.

For each $n$ the same argument can be applied to show that one of the affine coordinates of $f^n$ is a rational
function of degree at least $a \cdot c_+(N)^n$.  The same holds for the homogeneous expression for $f^n$, giving
$\da(f^n)  \geq a \cdot c_+(N)^n$, which is sufficient for the desired bound on $\lambda_1(f)$.
\end{proof}

\subsection{Equality on a dense set of full measure}\label{SUBSEC:EQUALITY_DENSE_OPEN_SET}

\begin{prop}\label{PROP:EQUALITY_DENSE_OPEN_SET}There is a dense set of full measure $\B'_N \subset \B_N$ so that $\lambda_1(f) = c_+(N)$ for $f \in \B'_N$.
\end{prop}

Before proving Proposition \ref{PROP:EQUALITY_DENSE_OPEN_SET}, we make some
comments about the action of \BPS on $\PP^2$.  It simplifies the discussion
to consider only \BPS $f$ for which the zeros are distinct and non-zero.  This
will be standing assumption in this subsection.

We begin by writing  $f$ in homogeneous coordinates $[Z:W:T]$, with $T=0$ corresponding
to the line at infinity with respect to the usual affine coordinates $(z,w)$.  We write 
\begin{eqnarray*}
f([Z:W:T]) = \left[f_1(Z,W,T):f_2(Z,W,T):f_3(Z,W,T)\right]
\end{eqnarray*}
\noindent
with
\begin{eqnarray}
\label{EQN_F_HOMOGENIOUS}
f_1(Z,W,T) &=& \theta_1 \prod_{i=1}^m (Z-a_iT) \prod_{i=1}^n (W-b_i T) \prod_{i=1}^p (T-Z \overline{c_i})
		\prod_{i=1}^q (T-W \overline{d_i})  \nonumber \\
f_2(Z,W,T) &=& \theta_2 \prod_{i=1}^p (Z-c_iT) \prod_{i=1}^q (W-d_i T) \prod_{i=1}^m (T-Z \overline{a_i})
                \prod_{i=1}^n (T-W \overline{b_i}) \\
f_3(Z,W,T) &=&  \prod_{i=1}^m (T-Z \overline{a_i}) \prod_{i=1}^n (T-W \overline{b_i}) \prod_{i=1}^p (T-Z \overline{c_i}) \prod_{i=1}^q (T-W \overline{d_i}).  \nonumber
\end{eqnarray}
\noindent
Since the zeros of $f$ distinct and non-zero,
no common factors occur in Equation (\ref{EQN_F_HOMOGENIOUS}).  Therefore,
$\da(f) = m+n+p+q$. 

We will call the lines $Z - a_i T=0$ the {\em zeros} of $A$
and denote the union of such lines by $Z(A)$.  Similarly, we will call the
lines $T-Z\overline{a_i}=0$  the {\em poles} of $A$, denoting the union of such
lines by $P(A)$.  The collections of lines $Z(B), P(B), Z(C), P(C), Z(D),$ and
$P(D)$ are all defined similarly.  Because of the standing assumption on the zeros of $f$, none of these lines coincide
with either the $z$ or $w$-axes.

From (\ref{EQN_F_HOMOGENIOUS}) we see that the ``vertical'' lines from $P(A)$ and the ``horizontal''
lines from $P(B)$ are collapsed to the point at infinity $[1:0:0]$.  In a similar way $P(C)$ and $P(D)$
are collapsed to $[0:1:0]$.

For the monomial map associated to $N$, the lines of zeros collapse to
$[0:0:1]$ and (often) the line at infinity $T=0$ collapses to either $[1:0:0]$
or $[0:1:0]$.  If the zeros of $f$ are distinct and non-zero, these lines no longer collapse.

The points of indeterminacy for $f$ are precisely the points for which all three coordinates of
(\ref{EQN_F_HOMOGENIOUS}) vanish.  In particular:

\begin{lem}If the zeros of $f$ are distinct and non-zero, then $f$ has $2$ points of indeterminacy on the line at
infinity: $[1:0:0]$ and $[0:1:0]$. 
\end{lem}

\begin{rmk}\label{RMK:NON_AS}Such $f$ are never algebraically stable on $\PP^2$:
As mentioned previously the lines in $P(A) \cup P(B)$ collapse under $f$ to
$[1:0:0]$ and the lines in $P(C) \cup P(D)$ collapse to $[0:1:0]$.
\end{rmk}

Each of the intersection
points from 
\begin{eqnarray*}
&&Z(A) \cap P(B),\,\,P(A) \cap Z(B),\\
&&Z(C) \cap P(D),\,\,P(C) \cap Z(D),\,\, \mbox{and}\\
&& (P(A) \cup P(B)) \cap (P(C) \cup P(D))
\end{eqnarray*}
\noindent 
is a point of indeterminacy.  There are $2(mn+pq)+(mq+np)$ such points of
indeterminacy in $\C^2$ and none of them lie on the $z$ or $w$-axes.

We write $I_f$ to denote the indeterminacy points of $f$ and
$C_f$ to denote the critical set of $f$ (within which are all of the collapsing
curves of $f$).  Let $P_f = P(A) \cup P(B) \cup P(C) \cup P(D)$ be the union of
all lines of poles for $f$.

\begin{lem}\label{LEM:ROTATION}
Given any $f$ and $g$ differing by rotations (but with the same zeros: $\sigma_1 = \sigma_2$),
we have the following:
\begin{itemize}
\item $I_{f} = I_{g}$, 
\item $C_{f} = C_{g}$, and
\item $P_{f} = P_{g}$.
\end{itemize}
\end{lem}

\begin{proof}
For each $f$ and $g$, the indeterminacy points 
are given by the points where the corresponding lift $F$ or $G$, respectively, to $\C^3$ has all three coordinates
vanishing.  The rotation multiplies the first two coordinates of each map by non-zero constants $\theta_1$
and $\theta_2$, hence has no affect on the indeterminacy points. 

Similarly, any rotation of $f$ by factors $\theta_1, \theta_2$ (non-zero) changes $\det(DF)$ by the non-zero
factor $\theta_1 \theta_2$, so the critical curves are unaffected.

The third item follows similarly.
\end{proof}


We now prove Proposition \ref{PROP:EQUALITY_DENSE_OPEN_SET}, defining  $\B'_N$ within the proof.

\begin{proof}[Proof of Proposition \ref{PROP:EQUALITY_DENSE_OPEN_SET}:]

The strategy of proof is as follows.  We begin by restricting that for any $f
\in \B'_N$, the zeros of $f$ are distinct and non-zero, so that the action of
$f$ on $\PP^2$ is as described above.  

We fix the zeros $\sigma$ (satisfying the above restriction) and let
$f_{(\theta_1,\theta_2)}$ be the mapping with zeros $\sigma$ and rotations
$(\theta_1,\theta_2)$.  According to Lemma \ref{LEM:ROTATION}, each of these
mappings will have the same indeterminacy set, critical set, and collection of
poles.  Let $I_f^0 \subset \C^2$ be the collection of finite indeterminate
points.

We will select a full-measure subset $\Omega \equiv \Omega_\sigma  \subset
\T^2$ so that if $(\theta_1,\theta_2) \in \Omega$, then any collapsing curve
(other than the lines of poles) does not have orbit landing in $I_f^0$ or on
one of the lines of poles.  (Landing on the lines of poles is dangerous since
they are mapped to the indeterminate points $[1:0:0]$ and $[0:1:0]$).  {\em It
is not clear that any such curves exist, but we are unable to rule them out in
the general case.}

We will then blow up  $\PP^2$ at $[1:0:0]$
and $[0:1:0]$ obtaining $\widetilde{\PP^2}$ and show that every such $f_{(\theta_1,\theta_2)}$ extends
to algebraically stable map on $\widetilde{\PP^2}$, allowing us to use the technique described in \S \ref{SEC:HOW_TO_COMPUTE_DYN_DEG}.

\vspace{0.05in}
Therefore, we let
\begin{eqnarray*}
\B'_N := \{f \in \B_N \,:\, \mbox{the zeros} \, \sigma \, \mbox{of} \, f \, \mbox{are distinct and non-zero, and} \, (\theta_1,\theta_2) \in \Omega_\sigma \}.
\end{eqnarray*}
We fix $\sigma$ (as in the definition of $\B'_N$) and construct $\Omega \equiv
\Omega_\sigma$.  For simplicity of exposition, we suppose that there is only
collapsing curve $C$ (other than the poles).
It can be generalized to a finite number of them in the
obvious way.  (Note that since $C$ is not a pole we have that $f(C) \in \C^2$.)

It is convenient to allow the rotations $\theta_1,\theta_2$ to be any complex
numbers, and later restrict that they each have modulus equal to one.   So long
as neither $\theta_1$ or $\theta_2$ is zero, the proof of Lemma
\ref{LEM:ROTATION} gives that the indeterminacy set and poles of $f$ remain
unchanged, allowing to denote them by $I_f^0$ and $P_f$, independent of
$\theta_1,\theta_2 \neq 0$.  

If either $\theta_1$ or $\theta_2$ is zero, then $f_{(\theta_1,\theta_2)}$
degenerates, and some poles and indeterminacy points may disappear.   However,
this degenerate map remains a  holomorphic map away from 
$I_f^0$ and $P_f$.

Let
\begin{eqnarray*}
\Psi_1 := \{(\theta_1,\theta_2) \in \C^2 \, : \, f_{(\theta_1,\theta_2)}(C) \not \in  (I_f^0 \cup P_f)  = \emptyset \}
\end{eqnarray*}
\noindent
and, inductively, let
\begin{eqnarray*}
\Psi_{n+1} := \{(\theta_1,\theta_2) \in \Psi_{n} \, : \, f^{n+1}_{(\theta_1,\theta_2)}(C)  \not \in  (I_f^0 \cup P_f)  = \emptyset \}.
\end{eqnarray*}
\noindent
We will show that each $\Psi_n$ is the complement of a (proper) analytic subset of $\C^2$ and that $(0,0) \in \Psi_n$.
The proof is by induction on $n$.

Notice that $\rho(\theta_1,\theta_2) = f_{(\theta_1,\theta_2)}(C)$ is a
holomorphic function defined on $\C^2$ with $\rho(0,0) = (0,0) \not \in (I_f^0
\cup P_f)$.  Therefore, the set of $(\theta_1,\theta_2)$ with
$\rho(\theta_1,\theta_2) \in  (I_f^0 \cup P_f)$ is a proper analytic subset of
$\C^2$.  We let $\Psi_1$ be its complement.

Suppose that $\Psi_n$ is the complement of a proper analytic set in $\C^2$ and
that $(0,0) \in \Psi_n$.  Let  $\varrho(\theta_1,\theta_2) =
f^{n+1}_{(\theta_1,\theta_2)}(C)$, which is holomorphic on $\Psi_n$.  Suppose
that $\varrho(\Psi_n) \cap (I_f^0 \cup P_f) \neq \emptyset$.  (Otherwise,
$\Psi_{n+1} = \Psi_n$ and we are done.) 

Since $\Psi_n$ is connected and $\varrho(0,0) = (0,0) \not \in
(I_f^0 \cup P_f)$ we see that $\varrho$ is a non-constant holomorphic function on
$\Psi_n$.   In particular, the set of $(\theta_1,\theta_2) \in \Psi_n$ having $\varrho(\theta_1,\theta_2) \in (I_f^0 \cup P_f)$
is a proper analytic subset of $\Psi_n$.  We let $\Psi_{n+1}$ be its complement.

We now let $\Omega_n = \Psi_n \cap \T^2$, which is the complement of a proper real-analytic subset.
Thus, the set $\Omega = \cap \Omega_n$ is the complement of a countable union of sets of measure zero,
and hence a set of total measure.  

\vspace{0.05in}
The set of zeros $\sigma$ that are distinct and non-zero is a dense set
of full measure in $\D^{m+n+p+q}$ and for each such $\sigma$, $\Omega_\sigma$
is of full measure in $\T^2$. It follows that $\B'_N \subset \B_N$ is also a
dense subset and, by Fubini's Theorem, of full measure.
\vspace{0.1in}

We now  blow up $[1:0:0]$ and $[0:1:0]$ obtaining $\widetilde{\PP^2}$.  Recall
that the blow up of $\mathbb{C}^2$ at $(0,0)$ is
\begin{eqnarray*}
\widetilde{\mathbb{C}}^2_{(0,0)} = \left\{((w,t),l) \in \mathbb{C}^2 \times \mathbb{P}^1 \mbox{ : } (w,t) \in l \right\}
\end{eqnarray*}
\noindent
There is a canonical projection $\pi:\widetilde{\mathbb{C}}^2_{(0,0)}
\rightarrow \mathbb{C}^2$ and the fiber $E_{(0,0)}= \pi^{-1}((0,0))$ is
referred to as the {\em exceptional divisor}.  See \cite{GH}.  In fact, this
definition is coordinate independent so that the notion of blowing up a complex
surface $X$ at a point $p \in X$ is well-defined.  The exceptional divisor
above $p$ will be denoted by $E_p$.

For $f \in \B'_N$ we check that $f$ extends continuously (and hence
holomorphically) to the blow-up at $[1:0:0]$.  The calculation at $[0:1:0]$ is
identical, and we omit it.  We write $f$ in the affine coordinates $w = W/Z$
and $t = T/Z$ so that the point of indeterminacy $[1:0:0]$ is at the origin
with respect to these coordinates.  

We work in the chart $(t,\lambda) \mapsto (\lambda t,t,\lambda) \in \widetilde{\mathbb{C}}^2_{(0,0)}$.  With domain
in this chart and codomain in the typical chart $(z,w) = (Z/T,W/T)$ we find that $f$ induces:
\begin{eqnarray*}
(t,\lambda) \mapsto 
\left(\prod_{i=1}^m \frac{1-a_i t}{t-\overline{a_i}}\prod_{i=1}^n \frac{\lambda t-b_i t}{t-\overline{b_i}\lambda t}\,\, , \,\,
	\prod_{i=1}^m \frac{1-c_i t}{t-\overline{c_i}}\prod_{i=1}^n \frac{\lambda t-d_i t}{t-\overline{d_i}\lambda t}\right)
\end{eqnarray*}
\noindent
so that the extension to $E_{[1:0:0]}$ is given by taking the limit $t \ra 0$:
\begin{eqnarray*}
\lambda \mapsto
\left(\prod_{i=1}^m \frac{-1}{\overline{a_i}}\prod_{i=1}^n \frac{\lambda-b_i}{1-\overline{b_i}\lambda} \,\, , \,\,
        \prod_{i=1}^p \frac{-1}{\overline{c_i}}\prod_{i=1}^q \frac{\lambda-d_i}{1-\overline{d_i}\lambda}\right).
\end{eqnarray*}
\noindent
\noindent
The calculation can also be done in the coordinates $\lambda' = \frac{1}{\lambda}$, where
one sees that the extension is continuous to all of $E_{[1:0:0]}$, hence
holomorphic.  (We are essentially using that none of the zeros $a_i$ or $c_i$
are equal to $0$.)  Since the extension is non-constant with respect to $\lambda$,
the extension of $f$ sends $E_{[1:0:0]}$ to a non-trivial rational curve.  

The blow-up at $[0:1:0]$ follows similarly and the extension of $f$ also sends
$E_{[0:1:0]}$ to a non-trivial rational curve.  We denote by
$\widetilde{f}:\widetilde{\PP^2} \ra \widetilde{\PP^2}$ this extension of $f$ to the space  $\widetilde{\PP^2}$ that is obtained
by doing both blow-ups.

\vspace{0.1in}
Having blown up $[1:0:0]$ and $[0:1:0]$ we will now observe that each of the lines
of poles from $P(A) \cup P(B)$ covers $E_{[1:0:0]}$ with non-zero degree and
each of the lines of poles from $P(C) \cup P(D)$ covers $E_{[0:1:0]}$ with non-zero
degree.  {\em In particular $\widetilde{f}$ does not collapse any of these lines to points.}  

The calculation is the same for each line, so we show it for $z = \frac{1}{\overline{a_1}}$.  If we parameterize this line by $w=W/T$, then the image in coordinate $\rho = \frac{w'}{t}$ (here $w' = W/Z$) can be found by substituting
$W=w,\, Z= \frac{1}{\overline{a_1}},$ and $T=1$ into the quotient $f_2(Z,W,T)/f_3(Z,W,T)$.  We obtain
\begin{eqnarray*}
\rho(w) = \frac{\prod_{i=1}^p(1/\overline{a_1}-c_i) \prod_{i=1}^q (w-d_i)}{\prod_{i=1}^p(1-\overline{c_i}/\overline{a_1})\prod_{i=1}^q(1-w\overline{d_i})}.
\end{eqnarray*} 
\noindent
This is a rational map of degree $q$ if $a_1 \neq c_i$ for all $i=1,\ldots,p$, which holds by
hypothesis that $f \in \B'_N$.  (In fact one can do the same calculation in the other coordinate charts
on the line $z = \frac{1}{\overline{a_1}}$ and on $E_{[1:0:0]}$, but the result will also be
a rational map of degree $q$ in those coordinates, as well.)

Similar calculations show that under $\widetilde{f}$, each of the lines of
poles from $P(A)$ covers $E_{[1:0:0]}$ with degree $q$ and each of the lines
the poles from $P(B)$ cover $E_{[1:0:0]}$ with degree $p$.  The poles from
$P(C)$ cover $E_{[0:1:0]}$ with degree $n$ and the poles from $P(D)$ with
degree $m$.  

\vspace{0.1in}

Let $\widetilde{X}$ be the blow-up of complex surface $X$ at point $p$ and
$\pi: \widetilde{X} \rightarrow X$ be the corresponding projection.  Given an
algebraic curve $D \subset X$ there are two natural ways to ``lift'' $D$ to
$\widetilde{X}$: the {\em total transform} and the {\em proper transform}.
The total transform is just $\pi^{-1}(D)$ while the proper transform is
obtained by the closure $\overline{\pi^{-1}(D \setminus \{p\})}$.  Clearly when
$p \not \in D$ there is no difference, however when $p \in D$ they differ by
the exceptional divisor $E_p \subset \widetilde{X}$.  In the case the many
points have been blown-up the analogous definitions hold,  see \cite{GH}.

\vspace{0.1in}

We now check that the only collapsing curves for
$\widetilde{f}:\widetilde{\PP^2} \ra \widetilde{\PP^2}$ are the proper
transforms of the curves $C_1,\ldots,C_k$ that are collapsed under $f$ to points in
$\C^2$. In fact any collapsing curve must be either the proper transform of a
collapsing curve for $f:\PP^2 \ra \PP^2$ or be one of the exceptional divisors
$E_{[1:0:0]}$ or $E_{[0:1:0]}$.  Since  $\widetilde{f}$ maps each of
$E_{[1:0:0]}$ or $E_{[0:1:0]}$ to a non-trivial rational curve, neither is a
collapsing curve.  Furthermore,  we have just checked that the lines from
$P(A)\cup P(B)\cup P(C) \cup P(D)$ are no longer collapsed by $\widetilde{f}$.
All that remains are the proper transforms of $C_1,\ldots,C_k$.

By the choice of $f \in \B'_N$ we have that the orbits of these
collapsing curves avoid the indeterminate points as well as all of the lines of
poles.  Therefore, under the extension $\widetilde{f}$, their orbits cannot
land on $E_{[1:0:0]}$, $E_{[0:1:0]}$, or the line at infinity.  Thus, the
orbits under $\widetilde{f}$ coincide with those under $f$, and they do not hit points in $I_f^0$, which are
the only indeterminate points for $\tilde{f}$.

\vspace{0.2in}
We can now
compute $\lambda_1(f) = \lambda_1(\widetilde{f})$ as the spectral radius of the
action of $\widetilde{f}^*$ on $H^{1,1}\left(\widetilde{\PP^2},\mathbb{R}\right)$.

\vspace{0.1in}
Let $\widetilde{L}_v \subset \widetilde{\PP^2}$ be the proper transform of the vertical
line $L_v := \{Z=0\}$ and let $\widetilde{L}_h$ be the proper transform of the horizontal line
$L_h := \{W=0\}$.  We choose the fundamental classes $[\widetilde{L}_v]$, $[E_{[0:1:0]}]$ and
$[E_{[1:0:0]}]$ as our basis of
$H^{1,1}\left(\widetilde{\PP^2},\mathbb{R}\right)$.  It will be useful in our
calculation to express $[\widetilde{L}_h]$ in terms of this basis.  

\begin{lem}\label{LEM:COH}
We have that:
\begin{eqnarray*}
[\widetilde{L}_h] \sim [\widetilde{L}_v] + [E_{[0:1:0]}] - [E_{[1:0:0]}].
\end{eqnarray*}
\end{lem}

\begin{proof}
Both $[L_v]$ and $[L_h]$ are cohomologous in $\PP^2$ so that
their total transforms $\pi^*([L_v]) = [\widetilde{L}_v] + [E_{[0:1:0]}]$ and $\pi^*([L_h]) = [\widetilde{L}_h] +
[E_{[1:0:0]}]$ are cohomologous, as well.
\end{proof}

We have that $f:\PP^2 \rightarrow \PP^2$ maps the lines of zeros $Z(A),Z(B)$
and the lines of poles $P(C), P(D)$ to $L_v$.  However, after blowing up
$[0:1:0]$ the lines of poles cover $E_{[0:1:0]}$ so that they should be
considered as part of $\widetilde{f}^*{[E_{[0:1:0]}]}$ and not part of
$\widetilde{f}^*{[\widetilde{L}_v]}$.

To see this more formally we write $\widetilde{f}$ with domain in the affine
coordinates $z=Z/T, w=W/T$ and image in the coordinates $t = T/W$, $\lambda =
Z/T$.  These image coordinates are chosen so that when $t=0$, $\lambda$ parameterizes $E_{[0:1:0]}$ (except for one point).
The second coordinate of the the image $(t,\lambda) =
\widetilde{f}(z,w)$ is given by:
\begin{eqnarray}
\lambda &=& \frac{\Pi_{i=1}^m (z-a_i ) \Pi_{i=1}^n (w-b_i)}{\Pi_{i=1}^m (1-z \bar a_i) \Pi_{i=1}^n (1-w \bar b_i)} \label{EQN:COMP_OF_PULLBACKS1}
\end{eqnarray}

In the $(t,\lambda)$ coordinates 
$\widetilde{L_v}$ is given by $\lambda=0$.   Therefore,
\begin{eqnarray}\label{EQN:F_LV}
\widetilde{f}^*(\widetilde{L_v}) \sim m[\widetilde{L_v}] + n
[\widetilde{L_h}]
\end{eqnarray}
\noindent
because of the $m$ factors of $z-a_i$ and $n$ factors of
$w-b_i$ in the numerator of Equation (\ref{EQN:COMP_OF_PULLBACKS1}).
Using Lemma \ref{LEM:COH} we find:
\begin{eqnarray*}\widetilde{f}^*{[\widetilde{L}_v]} \sim
(m+n)[\widetilde{L}_v] + n[E_{[0:1:0]}] - n[E_{[1:0:0]}].
\end{eqnarray*}

Suppose that 
we had parameterized the image of $\widetilde{f}$ in the other natural set of coordinates
in a neighborhood of $E_{[0:1:0]}$, given by $\hat z=
Z/W$ and $\eta = T/Z$.  Then, the total transform $\pi^{-1}(L_v) = 
\widetilde{L_v} \cup E_{[0:1:0]}$ is given by $\hat z=0$.  
The first coordinate of the image $(\hat z,\eta) = \widetilde{f}(z,t)$ is
given by
\begin{eqnarray*}
\hat z = \frac{ \theta_1 \prod_{i=1}^m (z-a_i) \prod_{i=1}^n (w-b_i)
\prod_{i=1}^p (1-z \overline{c_i}) \prod_{i=1}^q (1-w \overline{d_i})
}{\Pi_{i=1}^p (z-c_i) \Pi_{i=1}^q (w-d_i) \Pi_{i=1}^m (1-z \bar a_i)
\Pi_{i=1}^n (1-w \bar b_i)}
\end{eqnarray*}
\noindent
so that 
\begin{eqnarray*}
\widetilde{f}^* \left([L_v] +  [E_{[0:1:0]}]\right) \sim m[L_v] + n[L_h]+p[L_v]+q[L_h].
\end{eqnarray*}
\noindent
By substracting (\ref{EQN:F_LV}) and using Lemma \ref{LEM:COH} we find:
\begin{eqnarray*}
\widetilde{f}^*([E_{[0:1:0]}]) \sim (p+q)[\widetilde{L}_v]+q[E_{[0:1:0]}]-q[E_{[1:0:0]}].
\end{eqnarray*}

\noindent
A similar calculation gives  that
\begin{eqnarray*}
\widetilde{f}^*([E_{[1:0:0]}]) \sim (m+n)[\widetilde{L}_v] +n [E_{[0:1:0]}] - n [E_{[1:0:0]}].
\end{eqnarray*}

Therefore, in terms of the basis $\left\{[\widetilde{L}_v], [E_{[0:1:0]}],  [E_{[1:0:0]}] \right\}$ we have $\widetilde{f}^*$
given by:

\begin{eqnarray}\label{EQN:ACTION_H11}
\left[\begin{array}{ccc} (m+n) & (p+q) & (m+n) \\ n & q & n \\ -n & -q & -n \end{array} \right].
\end{eqnarray}

Therefore, $\lambda_1\left(\widetilde{f}\right) = r_1\left(\widetilde{f}^* \right)$ is the
largest eigenvalue of this matrix, which one can see coincides with $c_+(N)$.  Since dynamical degrees are invariant under birational conjugacy with $\widetilde{f}$ and $f$
conjugate under the projection $\pi$ we find $\lambda_1(f) = \lambda_1\left(\widetilde{f}\right) = c_+(N)$, as well.

This concludes the proof of Proposition \ref{PROP:EQUALITY_DENSE_OPEN_SET}.
\end{proof}

\begin{question}\label{QUESTION:STABLE_MODEL}
As mentioned earlier, in \cite{FAVRE_MON} it is shown that for any matrix of
degrees $N$ with positive coefficients there is some toric surface $\check{X}_N$ on
which monomial map corresponding to $N$ becomes algebraically stable.  Do all
$f \in \B_N$ extend to algebraically stable maps on $\check{X}_N$, as well?  This 
would be particularly helpful for studying bifurcations within the family.

In the case that $f$ is birational, \cite[Thm. 0.1]{DF_BIRATIONAL} gives the existence of
a modification by blow-ups $X$ of $\PP^2$ so that $f:X \rightarrow X$ is
algebraically stable.  Does $X = \widetilde{\PP^2}$ for birational \BPSNSP?  
\end{question}

\subsection{Equality everywhere}\label{SUBSEC:EQUALITY_EVERYWHERE}

\begin{lem}\label{LEM:GROWTH_OF_DEG}
Given $\epsilon > 0$ there exist $K_\epsilon$ so that for all $f \in \B'_N$ we have
\begin{eqnarray}\label{EQN:GROWTH_OF_DEG}
\da(f^n) \leq K_\epsilon (c_+(N) + \epsilon)^n.
\end{eqnarray}
\end{lem}

\begin{proof}
The action $\tilde f^* : H^{1,1}\left(\widetilde{\PP^2},\mathbb{R}\right)
\rightarrow H^{1,1}\left(\widetilde{\PP^2},\mathbb{R}\right)$ is given by the
matrix (\ref{EQN:ACTION_H11}) and hence independent of $f \in \B'_N$.
Therefore, the sequence of degrees $\{\da(f^n)\}$ is independent of $f
\in \B'_N$ and satisfies $\lim (\da(f^n))^{1/n} = c_+(N)$. This is
sufficient to give (\ref{EQN:GROWTH_OF_DEG}).
\end{proof}

\begin{proof}[Proof of Theorem \ref{THM:DYN_DEGREE}]
Given $\epsilon > 0$, let $K_\epsilon$ be given according to Lemma
\ref{LEM:GROWTH_OF_DEG}.  We now show that (\ref{EQN:GROWTH_OF_DEG}) actually
holds for every $f \in B_N$.  This will give 
$\lambda_1(f) \leq c_+(N)$. Combined with
the lower bound from Proposition \ref{PROP:DYN_DEG_LOWER}, it will complete the
proof of Theorem \ref{THM:DYN_DEGREE}.

Consider the family of $n$-th iterates $f^n$, where $f$ ranges over all of
$\B_N$.  If we write $f^n$ in homogeneous coordinates as $[h_1:h_2:h_3]$,
certain common factors will exist independent of the choice of 
$f \in \B_N$.  After eliminating all such common factors,
we obtain a homogeneous representation $[\hat h_1: \hat
h_2: \hat h_3]$ for $f^n$ that has no common factor for at least
one particular $f_0 \in \B_N$.

It is a consequence of elimination theory (see, e.g., \cite[\S 3.5]{IVA})
that $\hat h_1, \hat h_2,$ and $\hat h_3$ have a common factor if and only if
their coefficients satisfy an algebraic condition.  Since the coefficients of
the $\hat h_i$ are polynomial in the coefficients of $f = [f_1:f_2:f_3]$, 
the $\hat h_i$ have a common factor if and only if the
coefficients of $f$ satisfy an algebraic condition.

Since $f$ depends polynomially on $a_1,\ldots,d_q$, $\bar a_1,\ldots \bar d_q,
\theta_1,$ and $\theta_2$,  such a common factor occurs if and only if
$a_1,\ldots,d_q, \theta_1$, and $\theta_2$ satisfy a real-algebraic equation.
Because a common factor does not exist when representing $f_0^n$, 
this is a proper real-algebraic subset of $\B_N$.

By Proposition \ref{PROP:EQUALITY_DENSE_OPEN_SET}, $\B'_N$ is a dense 
in $\B_N$.  So, we find some $f_1 \in \B'_N$, so that $[\hat h_1:
\hat h_2: \hat h_3]$ represents of $f_1^n$ and has no common factor.  It
follows from Lemma \ref{LEM:GROWTH_OF_DEG} that by $\deg(\hat h_i) = \da(f_1^n)
\leq K_\epsilon (c_+(N) + \epsilon)^n.$

However, any $f \in \B_N$, has $[\hat h_1: \hat h_2: \hat h_3]$ as a homogeneous representation of $f^n$ (possibly with common factors). Therefore, $\da(f^n) \leq \deg(\hat h_i) \leq K_\epsilon (c_+(N) + \epsilon)^n.$
\end{proof}

\section{Case I: Large topological degree}
\label{SEC:LARGE_DEGREE}

As noted in the introduction, if $\det N > c_+(N)$, then every $f \in \B_N$
falls into Case I, having  $\dt(f) > \lambda_1(f)$.  We now check that for every
$N$, generically chosen \BPS are also from Case I.

Let $\hat \B_N$ be the set of \BPS for which all of the zeros from $\sigma$ are
distinct and none of the zeros are critical for their corresponding
one-variable Blaschke factor.  I.e. $A'(a_i) \neq 0$ for all $i$, and similarly
for $B,C$, and $D$.  It is straightforward that $\hat \B_N$ is an open dense
subset of $\B_N$ having total measure.  Furthermore it is invariant under
rotations by $\theta_1, \theta_2$.

\vspace{0.1in}
Recall:
\begin{thm*}{\bf \ref{THM:GENERICALLY_LARGE_DEGREE}}
For any matrix of degrees $N$
there is an open dense set of full measure $\hat \B_N \subset \B_N$ so that if
$f \in \hat \B_N$ then $\dt(f) = mq+np > \lambda_1(f)$.
\end{thm*}

\begin{proof}
It suffices to count the preimages of any point that is not a critical value of
$f$.  For $f \in \hat \B_N$, the origin $(0,0)$ is not a critical value.  This
follows because the preimages of $(0,0)$ in $\C^2$ are precisely the collection
of points from $Z(A) \cap Z(D)$ and $Z(B) \cap Z(C)$.  Substituting into $Jf =
A'(z)B(w)C(z)D'(w) - A(z)B'(w)C'(z)D(w)$ we see that the definition of $\hat \B_N$
prevents the $Jf$ from vanishing on these points.

Since the zeros are distinct for $f \in \hat \B_N$ we have $mq+np$ such points.  In
$\PP^2$ the line at infinity $T=0$ is forward invariant, so that there are no
additional preimages of $(0,0)$ that are not in $\C^2$.  Then, this total number of
preimages of the non-critical value $(0,0)$ is $\dt(f)$.

Notice that
\begin{eqnarray*}
\lambda_1(f) &=& \frac{m+q+\sqrt{(m-q)^2 +4np}}{2} <
\frac{m+q+\sqrt{(m-q)^2}+2\sqrt{np}}{2} \\ &=& \frac{m+q + |m-q|}{2} +
\sqrt{np} < mq+np = \dt(f).
\end{eqnarray*}
\end{proof}

For mappings in Case I of Conjecture \ref{CONJ:ENTROPY}, \cite[Thm
2.1]{GUEDJ_LARGE_DEGREE} gives a unique ergodic invariant measure $\mu$ of
maximal entropy $\log(\dt(f))$.  This measure is also backwards invariant,
satisfying $f^* \mu = \dt(f) \ \mu$.  (The notion of pulling back a measure is
not canonical; see \cite[p. 899]{RUSS_SHIFF} for the precise definition).  In particular,
$\supp(\mu)$ is totally invariant.

\begin{prop}\label{PROP:PREIMAGES_OF_T2}
There are points $x \in \T^2$ for which the weighted sequence of measures
\begin{eqnarray}\label{EQN:PREIMAGES}
\frac{1}{\left(\dt(f)\right)^n}{f^n}^* \delta_x
\end{eqnarray}
\noindent
converges weakly to $\mu$.  (Here $\delta_x$ indicates the Dirac mass.)
\end{prop}

\begin{proof}
It is a consequence of \cite[Thm 3.1]{GUEDJ_LARGE_DEGREE} that the sequence of
measures (\ref{EQN:PREIMAGES}) converges weakly to $\mu$, so long as $x$ is not
in a pluripolar exceptional set ${\mathcal E}_f$.  Since $\T^2$ is {\em
generating} (i.e. the complexification of each tangent space to $\T^2$ spans
the full tangent space in $\PP^2$), it follows from a theorem by Sadullaev
\cite[Theorem 4]{SAD} that ${\mathcal E}_f \cap \T^2$ has zero Haar measure in
$\T^2$.  
\end{proof}

If $\det N > c_+(N)$ then mappings $f \in \B_N$ having $\dt(f) = \det N$ can be
constructed by making an appropriate non-generic choice of zeros $\sigma$.

\begin{cor}\label{COR:MU_IN_T2}If $f \in \B_N$ satisfies $\dt(f) = \det N > c_+(N)$, then  $\supp(\mu) \subset \T^2$.
\end{cor}

\begin{proof}
According to Proposition \ref{PROP:PREIMAGES_OF_T2} the sequence of measures
(\ref{EQN:PREIMAGES}) converges weakly to $\mu$ when starting with a generic
point $x \in T^2$ (chosen with respect to the Haar measure).  Since $\dt(f) =
\det N = \dt(f_{|\T^2})$, all preimages of $x$ remain in $\T^2$ so that each of
the measures (\ref{EQN:PREIMAGES}) is supported in $\T^2$ and, therefore, $\mu$
is also.
\end{proof}

The hypothesis of Corollary \ref{COR:MU_IN_T2} are not satisfied for generic
\BPSNSP.  Rather, for any $f \in \hat \B_N$ we have $\dt(f) >
\lambda_1(f)$ and $\dt(f)
> \det N$.  In this case we have:

\begin{prop}\label{PROP:GENERICALLY_NO_CHARGE_T2}Suppose that $f \in \B_N$ satisfies $\dt(f) > \lambda_1(f)$ and $\dt(f) > \det N$, then $\mu$ does not charge $\T^2$.
In particular, this holds for any $f \in \hat \B_N$.
\end{prop}

\begin{proof}
Since $\dt(f) > \det N = \dt(f|\T^2)$ we have that $\T^2$ is not totally
invariant.  Since $\supp(\mu)$ is totally invariant, we cannot have $\supp(\mu)
\subset \T^2$.  However, since
$\T^2$ is forward invariant and $\mu$ is ergodic, we must have $\mu(\T^2) = 0$.
\end{proof}

Notice that 
even though $\mu$ does not charge $\T^2$,
its support may accumulate to $\T^2$.  In certain cases, we can rule out this possibility, using hyperbolic theory in a complex neighborhood of $\T^2$.

Suppose that the eigenvalues of $N$ satisfy $c_-(N) < 1 < c_+(N)$ so that
monomial map associated to $N$ induces a linear Anosov map with one-dimensional
stable direction and one dimensional unstable directions.  Then, any $f \in
\B_N$ with sufficiently small choice of zeros will also be an Anosov map of
$\T^2$ again with one dimensional stable and unstable directions.  Furthermore,
since $\T^2$ is hyperbolic for $f_{|\T^2}$ and $\T^2$ is generating, 
$\T^2$ is also a hyperbolic set for $f:\PP^2 \ra \PP^2$ with one-complex
dimensional stable and unstable directions.  For details on the hyperbolic
theory of endomorphisms, see \cite{JON_THESE}.

\begin{prop}\label{PROP:T2_ISOLATED}
Suppose that the eigenvalues of $N$ satisfy $c_-(N) < 1 < c_+(N)$ and that $f \in \B_N$.
If the zeros of $f$ are chosen sufficiently small, then $\T^2$ is isolated in the recurrent set of $f$.
\end{prop}

\begin{proof} 
We can assume that $\T^2$ is a hyperbolic set for $f:\PP^2 \ra \PP^2$.
Since $f_{|\T^2}$ is typically just an endomorphism, we go to the
natural extension $\hat \T^2 := \{(x_i)_{i \leq 0} \ : x_i \in \T^2 \, \mbox{and} \, f(x_i) = x_{i+1}\}$.  We denote
such histories by $\hat x = (x_i)_{i \leq 0} \in \hat \T^2$.
See \cite{JON_THESE}.  

Let us check that $\T^2$ is maximally invariant, i.e.
for a sufficiently small complex neighborhood $U$ of $\T^2$ we have:
\begin{eqnarray*}
\bigcap_{n \in \mathbb{Z}} f^n(U) = \T^2, \,\,
\bigcap_{n > 0} f^n(U) = W^u_{\rm loc}(\T^2), \,\, \mbox{and} \,\,
\bigcap_{n < 0} f^n(U) = W^s_{\rm loc}(\T^2).
\end{eqnarray*}
\noindent where 
$W^{s}_{\rm loc}(\T^2) = \cup_{x \in \T^2}
W^{s}_{\rm loc}(x)$ and $W^{u}_{\rm loc}(\T^2) = \cup_{\hat x \in \hat \T^2}
W^{u}_{\rm loc}(\hat x)$.

This is equivalent to the existence of a local
product structure for the natural extension $\hat \T^2$; see Definition 2.2 and
Corollary 2.6 from \cite{JON_THESE}.   Since $f_{|\T^2}$ is Anosov, there is
a local product structure (within $\T^2$) for $\hat \T^2$.  That is:
the unique point of intersection between $W^s_{\rm loc}(x)$ and $W^u_{\rm
loc}(\hat y)$ occurs at a point $z \in \T^2$ having some appropriate preorbit
$z_j$ (for $j < 0$) with $z_j \in W^u_{\loc}(\hat f ^{-j}(\hat q))$.
This local product structure naturally carries over when we consider $\T^2 \subset \PP^2$:
the unique point of intersection between the complex manifolds
$W^s_{\loc, \C}(x)$ and $W^u_{\loc,\C}(\hat y)$ must be the same point $z \in
\T^2$ and the previously chosen preorbit $z_j$ ($j < 0$) satisfies $z_j \in W^u_{\loc, \C}(\hat f
^{-j}(\hat q))$.

We suppose that there are recurrent points $r_i$ for $f:\PP^2 \ra \PP^2$ that
accumulate arbitrarily close to $\T^2$ from outside of $\T^2$.

In \cite[Prop. 3.8]{BP_PS} it was shown there are unique (semi) attracting
points $e \in \overline \D^2$ and $e' \in \overline{(\C \sm \D)^2}$ so that
$\D^2 \subset W^s(e)$ and $(\C \setminus \bar \D)^2 \subset W^s(e')$.  In
particular, these bidiscs contain no recurrent points other than $e$ and $e'$,
which are isolated from $\T^2$, since $f_{|\T^2}$ is Anosov.   Therefore, the
$r_i$ must accumulate to $\PP^2$ outside of $\D^2 \cup (\C \setminus \bar
\D)^2$ and their orbits cannot enter these bidiscs.  We will show that this is
impossible.

We can find some recurrent point $r_I$ in an arbitrarily small neighborhood $V
\subset U$ or $\T^2$.  Since $\T^2$ is maximally invariant and $r_I \not \in \T^2$, we must have $f^n(r_I) \not
\in U$ for some $n > 0$.  Let $n_0$ be the first such $n$.  If we choose $V$
sufficiently small, we can make $n_0$ arbitrarily large.
Therefore, we can assume that $f^{n_0
-1}(r_I)$ is arbitrarily close to  $W^u_{\loc, \C}(\hat x)$ for some $\hat x \in \hat \T^2$.

It is sufficient to check that $W^u_{\rm loc}(\hat x)  \subset \T^2 \cup \D^2
\cup (\C \setminus \bar \D)^2$.  Notice that $f^{n_0 -1}(r_I) \not \in
f^{-1}(U)$, which is some neighborhood of $\T^2$.  Thus, $f^{n_0-1}(r_I)$ is
arbitrarily close to $W^u_{\rm loc}(\hat x) \sm f^{-1}(U) \subset \D^2 \cup (\C
\setminus \bar \D)^2$, which will contradict the assumption that $r_I$ is
recurrent.

For all $x \in \T^2$, consider the complex conefield:
\begin{eqnarray*}
\KK(x) = \{v \in T_x \PP^2 \, : \, |\im(d\phi+d\psi)(v)| > |\im(d\phi - d\psi)(v)| \}.
\end{eqnarray*}
\noindent
which states precisely that $v$ is pointing into the pair of invariant bidiscs
$\D^2 \cup (\C \sm \D)^2$.  Therefore, forward invariance of $\KK$ follows from
forward invariance of these bidiscs.  Following general principles, since $\KK$
is invariant on the compact invariant set $\T^2$, we can extend it to an
invariant conefield in some small complex neighborhood of $\T^2$.

Because $\KK$ is forward invariant, the complex unstable manifold $W^u_{\rm
loc}(\hat x)$ of each $\hat x \in \hat \T^2$ is constrained within the cones.
This gives $W^u_{\rm loc}(\hat x) \subset \T^2 \cup \D^2 \cup (\C \sm
\D)^2$, as needed.
\end{proof}

\begin{rmk}
If $f_{|\T^2}$ is Anosov with two unstable directions (i.e. $1 < c_-(N) \leq c_+(N)$) then $\T^2$ is repelling.
In this case, points arbitrarily near to $\T^2$ can escape a neighborhood of $\T^2$ 
outside of $\D^2 \cup (\C \sm \D)^2$.  If $\dt(f) >
\dt(f_{\T^2})$, such an orbit could possibly approach a preimage of $\T^2$, allowing for
recurrence.  However, if $\dt(f) = \dt(f_{|\T^2})$, then elementary
considerations give that $\T^2$ is isolated in the recurrent set.
\end{rmk}

\begin{cor}Suppose that $f \in \B_N$ satisfies the hypothesis of Proposition \ref{PROP:T2_ISOLATED} and that $\dt(f) > \lambda_1(f)$, then
$\supp(\mu)$ is isolated from $\T^2$.
\end{cor}
\noindent
In particular, if $c_-(N) < 1 < c_+(N)$ and $f \in \hat \B_N$ with sufficiently small choice of zeros $\sigma$, then $\supp(\mu)$ is isolated from $\T^2$.

\begin{proof}Since $\dt(f) > \lambda_1(f) = c_+(N) > \det N$, Proposition \ref{PROP:GENERICALLY_NO_CHARGE_T2} gives that $\supp(\mu)$ is not contained in $\T^2$.  Therefore Proposition \ref{PROP:T2_ISOLATED} implies that they are isolated.
\end{proof}

\begin{cor}Suppose that the eigenvalues of $N$ satisfy $c_-(N) < 1 < c_+(N)$ and $f \in \hat \B_N$ has sufficiently small zeros $\sigma$. Then $\supp(\mu), \, \T^2,$ and $\{e,e'\}$
are isolated pieces of the recurrent set 
with $W^u(\supp(\mu)) \cap W^s(\T^2) \neq \emptyset$, $W^u(\T^2) \cap W^s(e) \neq \emptyset$, and  $W^u(\T^2) \cap W^s(e') \neq \emptyset$.
\end{cor}

\begin{proof}
By Proposition \ref{PROP:T2_ISOLATED}, $\supp(\mu)$ is isolated 
from $\T^2$ and in the proof we saw that $e$ and $e'$ are not in $\T^2$.
By Proposition
\ref{PROP:PREIMAGES_OF_T2}, preimages of $\T^2$ accumulate to $\supp(\mu)$ so that
$W^u(\supp(\mu)) \cap W^s(\T^2) \neq \emptyset$.  In the proof of Proposition
\ref{PROP:T2_ISOLATED} we saw that $W^u(\T^2)$ enters both $W^s(e)$ and $W^s(e')$.
\end{proof}

\begin{rmk}
Saddle sets for globally holomorphic maps were studied by \cite{DJ_SADDLE}, where it was shown
that the topological entropy of the saddle set is bounded above by
$\log(d_{\rm alg}(f))$, with equality holding if and only if the saddle set is
terminal.  The saddle set given by $\T^2$ that is discussed above conforms with
a possible generalization to meromorphic maps of the result of \cite{DJ_SADDLE},
since $\T^2$ is terminal and has topological entropy
$\log(\lambda_1(f))$ (because it is conjugate to the linear Anosov map).
\end{rmk}

\begin{question}
Suppose $N$ satisfy $c_-(N) < 1 < c_+(N)$, is there some open $U \subset \hat \B_N$ so that every $f \in U$ is Axiom-A, with non-wandering set $\Omega(f) = \supp(\mu) \cup \T^2 \cup \{e,e'\}$?
\end{question}

\section{Case II: Small topological degree}
\label{SEC:SMALL_DEGREE}

Given a choice of degrees $N$ with $c_+(N)> \det N$, the associated monomial
map has small topological degree (see Remark \ref{RMK:MONOMIALS}).
Non-monomial \BPS $f$ with $\lambda_1(f) > \dt(f)$ can also be constructed by
choosing the zeros $\sigma$ to have many repeated values.  We present one
specific family, for concreteness:

\begin{example}\label{EG:LOW_TOP_DEG}
For $a\neq b \neq c$, consider the family
\begin{eqnarray*}
f_{a,b,c}(z,w) = \left(\theta_1 \left(\frac{z-a}{1-\bar a z}\right)^5  \left(\frac{w-b}{1-\bar b w}\right)^2 ,
\,\, \theta_2  \frac{z-a}{1-\bar a z} \cdot \frac{z-c}{1-\bar c z} \cdot \frac{w-b}{1-\bar b w}\right).
\end{eqnarray*}

Members of this family are not in $\hat \B_N$ because $A(z)$ and $C(z)$ have a
common zero, however one can directly check that $\dt(f_{a,b,c}) = 5$.
Meanwhile, Theorem \ref{THM:DYN_DEGREE} gives  $\lambda_1(f) = c_+(N) =
\frac{6+\sqrt{32}}{2} > \dt(f_{a,b,c})$.

If $|a|, |b|,$ and $|c|$ are sufficiently small then $f_{a,b,c}$ is a
diffeomorphism on $\T^2$.   Therefore, Lemma \ref{PROP:MUTOR} gives an invariant
measure $\mutor$ supported on $\T^2$ with entropy $\log c_+(N) = \log
\lambda_1(f)$.  Since $\lambda_1(f) > \dt(f)$, \cite{DS_ENTROPY} gives that the
topological entropy of $f$ is $\log \lambda_1(f)$. Therefore, $\mutor$ is a
measure of maximal entropy for $f: \PP^2 \rightarrow \PP^2$.

Because $\dt(f) = 5$ there are many preimages of $\T^2$ in
$\PP^2$ that are away from $\T^2$.  This raises the question of whether there
exists a dynamically non-trivial invariant set outside of $\T^2$.
\end{example}

Recall 
\begin{prop*}{\bf \ref{PROP:MAX_ENTROPY_IN_T2}}
Let $f$ be a \BP diffeomorphism of small topological degree $\dt(f) < \lambda_1(f)$.  Then, $f:\PP^2 \rightarrow \PP^2$ has a measure of maximal entropy $\mutor$ supported within $\T^2$.
\end{prop*}

\begin{proof}
The general statement follows precisely as in Example \ref{EG:LOW_TOP_DEG}.
\end{proof}

\begin{question}
Is $\mutor$ is the {\em unique} invariant measure of maximal entropy  for $f:\PP^2 \ra \PP^2$?  
\end{question}

\begin{rmk}
It would be interesting to determine if the \BPS for which $\lambda_1(f) > \dt(f)$ satisfy the hypothesis of 
\cite{DDG1,DDG2,DDG3}.  
\end{rmk}

\begin{rmk}\label{REM:BIFURCATION} 

Suppose that $\det N < c_1(N)$.  Within $\B_N$, a change in the choice of zeros
$\sigma$ can result in the change between Case I and Case II, and therefore a
big change in the global dynamics on $\PP^2$.  For example, consider a
perturbation of the monomial map associated to $N$ by introducing arbitrarily
small generically chosen zeros.  Aside from degree considerations, it would be
interesting to know the mechanism(s) for this bifurcation.

This bifurcation may be similar to that obtained when applying a perturbation
of a monomial map, in a way that makes the resulting in a map that is globally
holomorphic on $\PP^2$.  Because the monomial map $f_0$ is hyperbolic on
$\T^2$, for small enough $\epsilon$ the perturbation to $f_\epsilon$ produces a
continuation of $\T^2$ to an invariant hyperbolic set for $f_\epsilon$.
However, the topological degree of $f_\epsilon$ jumps to $d^2$, so the
perturbation creates another invariant set of larger entropy $\log(d^2)$ (using
either \cite{GUEDJ_LARGE_DEGREE}, or earlier work as referenced in
\cite[\S 3]{S_PANORAME}).  

\end{rmk}

\section{Case III: Equal degrees}\label{SEC:EQUAL_DEGREES}

The monomial map associated to $N$ will have $\lambda_1(f) = 
\dt(f)$ if and only if $\det N =  c_+(N)$.  In this case we will have that
the smaller eigenvalue of $N$ satisfies $c_-(N) = 1$.  In appropriate coordinates, the linear
action of $f$ on $\T^2$ is given by product of an expanding map of $\T^1$ with the identity map of $\T^1$.

\begin{prop}\label{PROP:EQUAL_DEGREES}
Suppose that $f \in \B_N$ with $\lambda_1(f) = \dt(f)$, then  $c_+(N) = \det N = \dt(f)$ and $c_-(N) = 1$.
\end{prop}

\begin{proof}
Assume that there is some $f \in \B_N$ having
$\lambda_1(f) = \dt(f)$.  We write $c_\pm \equiv
c_\pm(N)$ for the eigenvalues of $N$.

Then,
\begin{eqnarray*}
c_- \cdot c_+ = \det N \leq \dt(f) = \lambda_1(f) = c_+
\end{eqnarray*}
\noindent
using $0 < \det N \leq \dt(f)$ and that $\lambda_1(f) = c_+$.
This gives $0 < c_- \leq 1$.

Yet, $c_-+c_+ = \tr N \in \N$ and $c_+ = \dt(f) \in \N$ so that $c_- \in \N$.
Hence $c_- = 1$ and $c_+ = \det N$.
\end{proof}

We can construct non-monomial \BPS $f \in \B_N$ with $\lambda_1(f) =
\dt(f)$ by choosing the zeros $\sigma$ to have many repeated values, as in \S
\ref{SEC:SMALL_DEGREE}.  

\begin{example}\label{EG:EQUAL_DEG}
Given zeros $a_1, a_2, a_3,$ and $b$, consider the family
{\small
\begin{align*}
& f_{a_1,a_2,a_3,b}(z,w) =  \\
& \,\,\,\,\,\, \left(\theta_1  \frac{z-a_1}{1-\bar a_1 z} \cdot \frac{z-a_2}{1-\bar a_2 z} \cdot \frac{z-a_3}{1-\bar a_3 z} \cdot \left(\frac{w-b}{1-\bar b w}\right)^2 ,
\,\, \theta_2  \frac{z-a_1}{1-\bar a_1 z} \cdot \frac{z-a_2}{1-\bar a_2 z}  \cdot \left(\frac{w-b}{1-\bar b w}\right)^3 \right)
\end{align*}
}
\noindent
corresponding to $N = \left[\begin{array}{cc} 3 & 2 \\ 2 & 3\end{array}\right]$.

\vspace{0.1in}
One can directly check
that $\dt(f_{a_1,a_2,a_3,b}) = 5 = \det N$.    Meanwhile,
Theorem \ref{THM:DYN_DEGREE} gives  $\lambda_1(f) = c_+(N) = 5$.
\end{example}

\begin{question}Suppose $f \in \B_N$ with $\lambda_1(f) = \dt(f)$.  Is the action of $f$ on $\T^2$ given by a skew product between an expanding map and a neutral map?
\end{question}

\appendix
\section{Entropy on $\T^2$}\label{APP:ENTROPY_T2}

\begin{prop}\label{PROP:ENTROPY_T2}
Let $f \in \B_N$ be a \BP diffeomorphism and $c_+(N)$ the leading eigenvalue of $N$. Then 
$h_{\rm top}(f_{|\T^2}) = \log(c_+(N)) > 0$.
\end{prop}

\begin{proof}
Central to the proof is that for any \BP
diffeomorphism $f$, the restriction $f_{|\T^2}$
has a global dominated splitting on all of $\T^2$;
see \cite[Cor 3.3]{BP_PS}. This places serious restrictions on the dynamics.
In particular, \cite[Thm 3.10]{BP_PS} gives that generic \BP diffeomorphisms
are Axiom-A with a very restricted behavior on the limit set.  If we
prove that $h_{\rm top}(f_{|\T^2}) = \log(c_+(N))$ for these maps, it will
follow for all \BP diffeomorphisms by the continuity of topological entropy for
$C^\infty$ surface diffeomorphisms \cite[Theorem 6]{Newhouse}.

Let $f \in \B_N$ be such a generic \BP diffeomorphism and let $f_0$ the
monomial map associated to $N$.  To simplify notation we will write $f \equiv
f_{|\T^2}$ and similarly for $f_0$.  Since $f_0$ is linear Anosov map induced
by $N$ we have $h_{\rm top}(f_0) = \det(c_+(N))$.

According to \cite[Lemmas 3.11 and  3.12]{BP_PS} there is a continuous
surjective $\pi:\T^2 \ra \T^2$ homotopic to the identity semiconjugating $f$ to
$f_0$ (and similarly for the lifts $\tilde{f}$ and $\tilde{f_0}$ to $\R^2$).  By
construction, we have that $\pi(x) = \pi(y)$ if and only if 
$\dist(\tilde{f}^n(x),\tilde{f}^n(y))$ remains bounded for all $n$.
Because of the semiconjugacy $\pi$, we have $h_{\rm top}(f) \geq h_{\rm top}(f_0)$.

Theorem 3.10 from \cite{BP_PS} gives that the limit set of $f$ consists of a
unique non-trivial homoclinic class $\HH$, and possibly a finite number
isolated of periodic points.  By ``unique homoclinic class'' we mean that,
given any saddle periodic point $p$, the closure of all all transverse
intersections of $W^s(p)$ with $W^u(p)$ is either $\HH$ or $\emptyset$.

The topological entropy of $f$ is concentrated on the limit set, so in this case
$h_{\rm top}(f) = h_{\rm top}(f_{|\mathcal{H}})$.  To complete the proof, we will check that $h_{\rm top}(f_{|\mathcal{H}}) \leq 
h_{\rm top}(f_{0})$.

Consider the restriction $\pi_{|\HH}:\HH \ra \pi(\HH)$ is a semiconjugacy onto its image.
Applying Bowen's Formula \cite[Theorem 17]{BOWEN}, we find
\begin{eqnarray}h_{\rm top}(f_{|\HH}) \leq h_{\rm top}(f_{0|\pi(\HH)}) + \sum_{z \in \pi(\HH)} h_{\rm top}\left(f_{|{\rm orbit \, of}\, \pi_{|\HH}^{-1}(z)}\right).
\end{eqnarray}
It suffices to check that for any $z \in \pi(\HH)$ we have 
$h_{\rm top}\left(f_{|{\rm orbit \, of}\, \pi_{|\HH}^{-1}(z)}\right) =0$.
To do this, we will show that if $x \neq y \in \pi_{|\HH}^{-1}(z)$ then both are in the
stable set of some periodic attracting interval $J$.  Note that
the entropy of a diffeomorphism of an interval is $0$.

Suppose that $\HH$ is generated by the periodic point $p$.  Since $\HH$ is
uniformly hyperbolic, there are local stable and unstable manifolds of a
uniform length over all of $\HH$.  Furthermore, for any $h \in \HH$, global
manifolds can be formed, e.g. $W^s(h) = \cup f^{-n}(W^s_{\rm loc}(f^n(h))$.
Notice that $W^u(p)$ accumulates to $h$, intersecting 
$W^s(h)$ transversely on at least one side $W^{s+}(h_1)$.  It follows from the
$\lambda$-Lemma \cite{PDM} that $W^{s+}(h)$ is unbounded and accumulates to every other
point $h_2 \in \HH$.  In particular, for any pair of points $h_1 \neq h_2$ in
$\HH$ there is a compact connected arc $I \subset W^{s+}(h_1)$ intersecting
$W^u_{\rm loc}(h_2)$.  Because 
$W^u(p)$ intersects $I$, the $\lambda$-lemma implies that $f^{-n}(I)$ has length arbitrarily large.  Further, since $\tilde f$ inherits a dominated splitting on
$\R^2$, this implies that the lifts $\tilde{f}^{-n}(I)$ will have arbitrarily
large diameter.

We apply the above discussion to the pair $x \neq y \in \pi_{|\HH}^{-1}(z)$, letting $x' \in
W^{s+}(x) \cap W^u_{\rm loc}(y)$, $I \subset W^{s+}(x)$ the arc connecting $x$
to $x'$, and $I' \subset W^u_{\rm loc}(y)$ be the arc connecting $x'$ to $y$.
Notice that $x' \neq y$ because the discussion from the previous paragraph would give
that
$\tilde{f}^{-n}(x)$ and $\tilde{f}^{-n}(y)$ become arbitrarily far apart in $\R^2$.
Similarly, that $f^n(I')$ remains a finite length because $x' \in
W^s(x)$ and $\tilde{f}^{n}(x)$ and $\tilde{f}^{n}(y)$ remain at finite distance
in $\R^2$.

We can now apply the Denjoy Property from \cite[\S 2.4]{PS2} to the interval
$I'$.  Since $I'$ is part of $W^u_{\rm loc}(y)$, it is tangent to the
center-unstable linefield $F$ from the dominated splitting.  In the
language of \cite[\S 2.4]{PS2}, $I'$ is a $\delta$-$E$-arc, where $\delta$ is the bound on the length of $f^n(I')$.  
Since there is a global dominated splitting for $f$
on all of $\T^2$, Theorem 2.3 from \cite{PS2} applies for $\delta$-$E$-arcs for any $\delta > 0$.  This gives that
$\omega(I')$ is either a periodic closed curve, a periodic closed arc $J$ (with $I' \subset W^s(J)$), or a periodic point which is either a sink or a saddle-node.  
By \cite[Thm 3.10]{BP_PS} there can be no periodic closed curves under $f$.  A sink is impossible because $x,y \in \HH$ cannot be in the basin of attraction of a sink.
A saddle-node is impossible because $f$ is Axiom-A.  Therefore, $I' \subset W^s(J)$ for some periodic closed arc $J$.
In particular $x',y \in I' \subset W^s(J)$ and $x' \in W^s(x)$ giving $x \in W^s(J)$.
\end{proof}

\begin{rmk}
If the set of all \BP diffeomorphisms within $\B_N$ is connected, a much simpler proof of Proposition \ref{PROP:ENTROPY_T2} would follow directly from \cite[Thm. E]{PS2}.
However, this is presently unknown.
\end{rmk}

\begin{prop}\label{PROP:MUTOR}
Let $f \in \B_N$ be a \BP diffeomorphism. Then, there is a unique measure
of maximal entropy $\mutor$ for $f_{|\T^2}$.
\end{prop}

\noindent
Note that here the meaning of ``maximal'' is with respect to invariant measures
supported on $\T^2$.  We have already observed in Proposition \ref{PROP:GENERICALLY_NO_CHARGE_T2} that the generic \BPS $f \in \hat B_N$ have
an invariant measure $\mu$ of higher entropy $\log(\dt(f)) > \log(c_+(N))$ that does not charge $\T^2$.

\begin{proof}
From \cite[Thm 3.10]{BP_PS} it follows that generic \BP diffeomorphisms are
Axiom-A on $\mathbb{T}^2$ with a unique non-trivial homoclinic class.  From
\cite[Thm E]{PS2}, the restriction of any \BP diffeomorphism to its limit set 
is conjugate to the restriction of one of these Axiom-A maps to its limit set.  
Therefore, on one hand, we conclude that
any  Blaschke product diffeomorphism has a unique non-trivial homoclinic class in the torus so that the topological entropy in $\mathbb{T}^2$ (which by
Proposition \ref{PROP:ENTROPY_T2} is $\log(c_+(N ))$) is equal to the topological entropy of
the diffeomorphisms restricted to this homoclinic class.  On the other hand,
since the homoclinic class is conjugate to a hyperbolic one, it follows that it
has a unique measure of maximal entropy with support in the homoclinic class.
\end{proof}

\bibliographystyle{plain}
\bibliography{BP.bib}

\end{document}